\documentclass{article}
\usepackage{amsmath}
\usepackage{amssymb}
\usepackage[english] {babel}

\def\a{\alpha} \def\b{\beta} \def\d{\delta} \def\f{\varphi}
\def\l{\lambda}  \def\R{\mathbb{R}} \def\C{\mathbb{C}}
\def\H{\mathbb{H}} \def\O{\mathbb{O}} \def\I{\mathbb{I}}
 \def\minus{\smallsetminus}
\def\({\left(} \def\){\right)} 
\def\<{\langle} \def\>{\rangle}
\def\w{\wedge} \def\bx{\hspace*{\fill}$\Box$\vspace{1ex}} 
 \def\Im{\mathrm{Im}}
\def\Aut{\mathrm{Aut}} 
\def\G{\mathcal{G}_{2}}
\def\pd{\mathfrak{p}(\delta)} \def\E{\mathcal{E}}
\def\p{\mathfrak{p}} \def\sp{\mathrm{span}}

\newcommand{\ie}{i.e.}

\newtheorem{thm}{Theorem}[section]
\newtheorem{lma}[thm]{Lemma}
\newtheorem{prop}[thm]{Proposition}

\title{Normal forms for the $\G$-action on the real symmetric $7\times 7$-matrices by
  conjugation}
\author{Erik Darp\"{o}\\ {\small \textit{Matematiska Institutionen, Uppsala Universitet,}}\\
  {\small\textit{Box 480, S-75106 Uppsala, Sweden.}}\\{\small\texttt{erik.darpo@math.uu.se}}}

\begin{document}
\selectlanguage{english}
\date{} 
\maketitle 

\begin{abstract}
The exceptional Lie group $\G\subset \mathrm{O}_7(\R)$ acts on the set of real symmetric
$7\times7$-matrices by conjugation. We solve the normal form problem for this group action. In
view of the earlier results \cite{bbo}, \cite{malaga} and \cite{coll}, this gives rise to a
classification of all finite-dimensional real flexible division algebras. By a
classification is meant a list of pairwise non-isomorphic algebras, exhausting all
isomorphism classes.

We also give a parametrisation of the set of all real symmetric matrices, based on eigenvalues.
\end{abstract}

\noindent
Mathematics Subject Classification 2000: 15A21, 17A20, 17A35, 17A36, 17A45. \\
Keywords: Normal form, group action, vector product, octonion, automorphism, real division algebra,
flexible algebra.

\section{Introduction} \label{intro}

A \emph{vector product} on a Euclidean space $V=(V,\<\,\>)$ is a linear
map $\pi:V\w V\to V$ with the property that the set $\{u,v,\pi(u\w v)\}\subset V$ is orthonormal
whenever $\{u,v\}\subset V$ is. 
A morphism $\pi\to\pi'$ of vector products $\pi$ and $\pi'$ on $V$ and $V'$ respectively
is an algebra morphism $(V,\pi)\to(V',\pi')$
preserving the scalar product. Vector products can be defined on Euclidean spaces of
dimension 0, 1, 3 and 7 only, and are unique up to isomorphism in each dimension.

If $\pi$ is a vector product on a 7-dimensional Euclidean space $V$, then $\Aut(\pi)$ is
isomorphic to the exceptional Lie group $\G$. We view this isomorphism as an
identification, and write $\G=\Aut(\pi)$. Since $\G$ is a subgroup of $\mathrm{O}(V)$,
\begin{equation}\label{nform}
  \mathrm{Sym}(V)\times\G\to \mathrm{Sym}(V),\;(\d,g)\mapsto g^{-1}\d g
\end{equation}
defines an action of $\G$ on the set 
$\mathrm{Sym}(V)$ of symmetric linear endomorphisms of $V$. The present article is devoted to
the solution of the normal form problem for this group action.
Our study is motivated by the classification theory for flexible division algebras, in
which this group action plays an important role (see Proposition~\ref{cuenca}).

A (not necessarily associative) algebra $A$ is said to be a \emph{division algebra} if
$A\neq\{0\}$ and the linear endomorphisms $\mathrm{L}_a:x\mapsto ax$ and $\mathrm{R}_a:x\mapsto xa$ are
bi\-ject\-ive for all $a\in A\minus\{0\}$. It is called \emph{alternative} if any
subalgebra generated by two elements is associative, \emph{power associative} if any
subalgebra generated by one element is associative, and \emph{flexible} if the identity
$x(yx)=(xy)x$ holds for all $x,y\in A$.
In this article, our attention is restricted to finite-dimensional real algebras
(henceforth referred to simply as 'algebras'). In the finite-dimensional case, the
division property is equivalent to $xy=0\Rightarrow x=0 \mbox{ or } y=0$ for all 
$x,y\in A$. 

The most well-known division algebras are the real numbers $\R$, the complex numbers
$\C$, the quaternion algebra $\H$ (Hamilton 1843) and the octonion algebra $\O$ (Graves
1843, Cayley 1845).
Classical theorems assert that $\{\R,\C,\H\}$ and $\{\R,\C,\H,\O\}$ classify
all associative and alternative division algebras respectively (Frobenius
\cite{frobenius}, Zorn \cite{zorn}), and that every division algebra has dimension either
1, 2, 4 or 8 (Bott and Milnor \cite{bottmilnor}, Kervaire \cite{kervaire}).

An algebra $A$ is said to be \emph{quadratic} if it has an identity element $1\neq 0$,
and the set $\{1,x,x^2\}$ is linearly dependent for all $x\in A$. It is known that a
real division algebra is quadratic if and only if it is power associative
(this is a consequence of \cite[Lemma~5.3]{zur}). Hence, in particular every
alternative division algebra is quadratic.

In any quadratic algebra $B$, the subset
\[ \Im B=\{ b\in B\minus\R1 \mid b^2\in\R1 \} \cup \{0\}\subset B \]
of purely imaginary elements is a linear subspace of $B$, and $B=\R1\oplus \Im B$
(Frobenius \cite{isomorphie}). We shall write $\alpha+v$ instead of $\alpha 1+v$ when
referring to elements in this decomposition.

A linear map $\eta:V\w V\rightarrow V$, where $V$ is a finite-dimensional Euclidean space,
is called a \emph{dissident map} on $V$ if the set $\{v,w,\eta(v\w w)\}\subset V$ is linearly
independent whenever $\{v,w\}\subset V$ is. If in addition $\xi:V\w V\rightarrow \R$ is 
a linear form, $(V,\xi,\eta)$ is called a \emph{dissident triple}. The class of dissident
maps is given the structure of a category, denoted $\mathcal{D}$, by declaring as
morphisms $(V,\xi,\eta)\to(V',\xi',\eta')$ those linear maps $\sigma:V\to V'$ for which
$\sigma\eta=\eta'(\sigma\wedge\sigma)$, $\xi=\xi'(\sigma\w\sigma)$ and  
$\<x,y\>=\<\sigma(x),\sigma(y)\>$ for all $x,y\in V$. The assignment $\pi\mapsto(V,0,\pi)$
defines a full embedding of the category of vector products into\nolinebreak[2] $\mathcal{D}$.

Each dissident triple $(V,\xi,\eta)\in\mathcal{D}$ determines a quadratic division algebra
$\mathcal{H}(V,\xi,\eta)= \R\times V$ with multiplication
\begin{equation*}
  (\a,v)(\b,w)=(\a\b-\<v,w\>+\xi(v\w w),\a w+\b v+\eta(v\wedge w)).
\end{equation*}
On defining $(\mathcal{H}\sigma)(\a,v)=(\a,\sigma(v))$ for morphisms
$\sigma:(V,\xi,\eta)\rightarrow(V,\xi',\eta')$, $\mathcal{H}$ becomes a functor
from $\mathcal{D}$ to the category $\mathcal{Q}$ of quadratic division
algebras (morphisms in $\mathcal{Q}$ are algebra morphisms preserving the identity
element).
This functor turns out to be an equivalence of categories (Osborn \cite{osborn},
cf. Dieterich \cite{dissalg}). Flexible quadratic division algebras correspond under
$\mathcal{H}$ to triples of the form $(V,0,\eta)$ for which $\<v,\eta(v\w w)\>=0$ for all
$v,w\in V$ (or, equivalently, $\<\eta(u\w v),w\>=\<u,\eta(v\w w)\>$ for all $u,v,w\in V$).

The categories of flexible quadratic division algebras and the corresponding dissident
maps are denoted by $\mathcal{Q}^{fl}$ and $\mathcal{D}^{fl}$ respectively. We write
only $\eta$ as shorthand for $(V,0,\eta)\in\mathcal{D}^{fl}$. Vector products correspond
to alternative division algebras under $\mathcal{H}$.

In \cite{bbo} and \cite{coll}, the classification problem for the flexible division
algebras is reduced to the classification problem for the subclass consisting of all
8-dimensional quadratic flexible division algebras. The latter problem is
addressed by Cuenca Mira et al.\ in \cite{malaga}. Our Proposition~\ref{cuenca} states
their main theorem in the language of dissident maps. Here $\mathrm{Pds}(V)$ denotes the
set of positive definite symmetric endomorphisms of the Euclidean space $V$, and
$\epsilon^\ast$ the adjoint of the linear endomorphism $\epsilon$.

\begin{prop}\label{cuenca} \cite[p. 21]{malaga}
  Let $\pi:\R^7\w\R^7\rightarrow\R^7$ be a vector product, and $\eta$ a flexible
  dissident map on a Euclidean space $V$ of dimension 7. Then the following holds.
  \begin{enumerate}
  \item For any $\epsilon\in \mathrm{GL}(V)$, $\epsilon^\ast\eta(\epsilon\w\epsilon)$ is a
    flexible dissident map.
  \item $\eta\cong\d^\ast\pi(\d\w\d)=\d\pi(\d\w\d)$ for some $\d\in \mathrm{Pds}(\R^7)$.
  \item For $\d_1,\d_2\in \mathrm{Pds}(\R^7)$,
    $\d_1\pi(\d_1\w\d_1)\cong\d_2\pi(\d_2\w\d_2)$ if and only if
    $\d_1=\sigma^{-1}\d_2\sigma$ for some $\sigma\in\Aut(\pi)$.
  \end{enumerate}
\end{prop}

This result reduces the problem of classifying the 8-dimensional flexible quadratic
division algebras to the normal form problem for the group action
\begin{equation}\label{pdsnform}
  \mathrm{Pds}(\R^7)\times\G\rightarrow \mathrm{Pds}(\R^7), \; (\d,g)\mapsto
  \d\cdot g=g^{-1}\d g.
\end{equation}
This is a subproblem of the normal form problem for the group action (\ref{nform}).

In Section~\ref{normal}, Propositions~\ref{p7}--\ref{p111duo}, normal forms for
(\ref{nform}) are given separately for each possible con\-fig\-ur\-ation of
eigenvalues. The solution for the positive definite case is obtained simply by restricting
attention to positive eigenvalues. In Section~\ref{parametrisering} we give
parametrisations of the sets of all symmetric matrices sharing some fixed con\-fig\-ur\-ation
of eigenvalues, some of which are needed in Section~\ref{normal}.

We use the following notation and conventions.
Given $m\in\mathbb{N}$, we set
$\underline{m}=\{\mu\in\mathbb{N} \mid 1\leqslant\mu\leqslant m\}$.
The identity matrix of size $n\times n$ is denoted by $\I_n$, and the identity map on a
set $X$ by $\I_X$.

The vector space $\R^n$ is viewed as an Euclidean space, equipped with the standard
scalar product. Elements in $\R^n$ are considered as column vectors, and linear
endomorphisms of $\R^n$ are identified with $n\times n$-matrices in the natural way.
The transpose of a matrix (or column vector) $A$ is denoted by $A^\ast$. 
We write $\underline{e}=(e_1,\ldots,e_n)$ for the standard basis in $\R^n$, \ie,
$e_k=(\d_{kl})_{l\in\underline{n}}$. 
If $\l$ is an eigenvalue of a linear endomorphism $\d$ of some vector space $V$, then
$\E_\l(\d)\subset V$ denotes the corresponding eigenspace.
For $u\in V\minus\{0\}$, where $V$ is a Euclidean space, we define
$\hat{u}=u^\wedge=\frac{1}{\|u\|}u$.
Given any subset $W\subset V$, $W^\perp$ denotes its orthogonal complement in $V$.
If $W\subset V$ is a subspace, then $P_W:V\rightarrow W$ denotes the projection onto $W$
along $W^\perp$. For $v\in V$, we write $v^\perp=\{v\}^\perp$ and
\begin{equation*}
  \<W,v\>=\max\{\<w,v\>\mid w\in W,\:\|w\|=1\}=\left\|P_W(v)\right\|.  
\end{equation*}

\section{The vector product in 7-dimensional space}\label{vektorprodukten}

In the following two sections, $V$ denotes a fixed 7-dimensional Euclidean
space, equipped with a vector product $\pi$. The map $\pi$ is considered as a
multiplication on $V$, and we write $xy=\pi(x\w y)$. Accordingly, $\mathrm{L}_x$ and $\mathrm{R}_x$
denote the maps $v\mapsto\pi(x\w v)$ and $v\mapsto\pi(v\w x)$ respectively.

The following lemma provides means to control the multiplication in the algebra $(V,\pi)$.

\begin{lma}\label{vm} Let $u,v\in V$ be orthonormal vectors. The following identities
  hold:
  \begin{enumerate}
  \item $u(uv)=-v$,
  \item $v(uv)=u$.
  \end{enumerate}
  In particular, $\pi$ induces a vector product on $\sp\{u,v,uv\}$. \\
  If in addition $z\in V$ is a unit vector orthogonal to $u$ and $v$, then
  \begin{enumerate}\setcounter{enumi}{2}
  \item $u(vz)=-(uv)z=(vu)z$.
  \end{enumerate}
\end{lma}

\textit{Proof:} If $x\in u^\perp$, then $\|ux\|=\|x\|$. This means that the
linear map $\mathrm{L}_u:u^\perp\rightarrow u^\perp$ is an isometry. Thus,
\begin{equation*}
  \<x,u(uv)\>=\<xu,uv\>=-\<ux,uv\>=\<ux,u(-v)\>=\<x,-v\>
\end{equation*}
for all $x\in V$, and hence $u(uv)=-v$. The second identity follows from the first via
anti-com\-muta\-tiv\-ity of $\pi$.

By \textit{1}, we also have $u(uz)=-z$. Moreover,
\begin{align*}
  -2v&=-\|u+z\|^{2}v=(u+z)\((u+z)v\)=(u+z)(uv+zv) \\
  &=u(uv)+u(zv)+z(uv)+z(zv)\\
  &=-v+u(zv)+z(uv)-v=-2v-u(vz)-(uv)z.
\end{align*}
Hence $u(vz)=-(uv)z$. \bx

Let $S\subset V$ be a $3$-dimensional oriented subspace.
We define $m:S\rightarrow V$ by $m(0)=0$ and $m(v)=\|v\|ef$ if $v\neq0$, where $e,f\in S$
are such that $(e,f,\hat{v})$ is a positively oriented orthonormal basis for $S$.
Given $v\in S\minus\{0\}$, for any two equally oriented orthonormal bases $(e,f)$ and
$(e',f')$ for $S\cap v^\perp$ we have $ef=e'f'$. Therefore, the map $m$ is well
defined.

The following lemma is an important tool for the solution of the normal form problem for
the group action (\ref{nform}).
\begin{lma}\label{m}
  \begin{enumerate}
  \item The map $m:S\to V$ is linear and orthogonal.\label{m1}
  \item The functions $S\to\R_{\geqslant0},\;v\mapsto\<S,m(v)\>$ and
  $S\to\R_{\geqslant0},\;v\mapsto\<S^\perp,m(v)\>$ are constant on the unit sphere in
  $S$.\label{m2}
  \item The map $S\to V,\;v\mapsto v(m(v))$ is constant on the unit sphere in
  $S$. \label{m3}
  \end{enumerate}
\end{lma}

\textit{Proof:} \ref{m1}. Let $v\in S\minus\{0\},\;\a\in\R\minus\{0\}$. Take $e,f\in S$ as
in the definition of $m$. Now
\begin{alignat*}{2}
  m(\a v)&=\|\a v\|ef=|\a|\|v\|ef=\a m(v) &&\mbox{ if } \a>0\\
  \intertext{and also}
  m(\a v)&=\|\a v\|fe=-\|\a v\|ef=-|\a|\|v\|ef=\a m(v) &&\mbox{ if } \a<0
\end{alignat*}
and certainly $m(0v)=m(0)=0=0m(v)$.
Given $u,v\in S\minus\{0\}$, we take $e\in S\cap u^\perp\cap v^\perp$, and $f,f'\in S$ such
that $(e,f,\hat{u})$ and $(e,f',\hat{v})$ are positively oriented orthonormal bases for
$S$. Thus $\hat{v}=a\hat{u}+bf$ and $f'=-b\hat{u}+af$ for some $a,b\in\R$ satisfying
$a^2+b^2=1$.
Now set $g=\|u\|f+\|v\|f'$. We have 
\begin{align*}
  \<u+v,g\>&=\|v\|\<u,f'\>+\|u\|\<v,f\>\\
  &=\|u\|\|v\|(\<\hat{u},f'\>+\<\hat{v},f\>)=\|u\|\|v\|(-b+b)=0\\
  \intertext{and}
  \|g\|&=\|u+v\|.\\
  \intertext{Using this, we calculate}
  m(u+v)&=\|u+v\|e\hat{g}=eg=\|u\|ef+\|v\|ef'=m(u)+m(v).
\end{align*}
This proves that $m$ is linear. Since $m$ maps unit vectors to unit vectors, it now follows
that it is also orthogonal.

\ref{m2}. Here it suffices to show that $v\mapsto\<S,m(v)\>$ is constant, since for all
$v\in S$ we have $\|v\|^2=\|m(v)\|^2=\<S,m(v)\>^2+\<S^\perp,m(v)\>^2$.
Let $\underline{b}=(b_1,b_2,b_3)$ be a positively oriented orthonormal basis for $S$.
Now $m(b_1)=b_2b_3\in\{b_2,b_3\}^\perp$, from which follows that
$P_Sm(b_1)=\<b_1,b_2b_3\>b_1$. Similarly,
$P_Sm(b_2)=\<b_2,b_3b_1\>b_2$ and $P_Sm(b_3)=\<b_3,b_1b_2\>b_3$.
Since 
\begin{equation*}
\<b_1,b_2b_3\>=\<b_3,b_1b_2\>=\<b_2,b_3b_1\>, 
\end{equation*}
for an arbitrary unit vector $v=v_1b_1+v_2b_2+v_3b_3$ we have
\begin{multline*}
\<S,v\>=\|P_Sm(v)\|=\|v_1P_Sm(b_1)+v_2P_Sm(b_2)+v_3P_Sm(b_3)\| \\
=|\<b_1,b_2b_3\>|\|v_1b_1+v_2b_2+v_3b_3\|=|\<b_1,b_2b_3\>|\|v\|=|\<b_1,b_2b_3\>|.
\end{multline*}
So the map $v\mapsto\<S,v\>$ is constant on the unit sphere of $S$.

\ref{m3}. If $S$ is closed under the multiplication map $\pi$, then $m(v)=\pm v$ and
$vm(v)=0$ for all $v\in V$, so our 
assertion holds true. Suppose $S$ is not closed under $\pi$.
Then $\<S^\perp,m(x)\>\neq0$ for all $x\in S\minus\{0\}$.
For any unit vector $u\in S$, add $v$ and $w$ such that $(u,v,w)$ becomes positively
oriented and orthonormal in $S$.
Now $\<v,um(u)\>=\<v,u(vw)\>=\<vu,vw\>=\<u,w\>=0$ and 
$\<m(v),um(u)\>=\<wu,u(vw)\>=\<(wu)u,vw\>=\<-w,vw\>=0$.
Analogously, we get $\<w,um(u)\>=0$ and $\<m(w),um(u)\>=0$.
Since also $\<u,um(u)\>=\<m(u),um(u)\>=0$, this shows that
$um(u)\in\{u,v,w,m(u),m(v),m(w)\}^\perp=(S+m(S))^\perp$.

From \ref{m1} and \ref{m2} follows that $m(S)\subset V$ is a 3-dimensional subspace,
and $S\cap m(S)=\{0\}$. Hence $\dim{(S+m(S))^\perp}=1$.
We also have $\|um(u)\|=\<S^\perp,m(u)\>$. Therefore, since $u\mapsto\<S^\perp,m(u)\>$ is
constant, there exists a $z\in V$, with $\|z\|=\<S^\perp,m(u)\>\neq0$, such that $um(u)=\pm z$
for all unit vectors $u\in S$. As the map $u\mapsto um(u)$ is continuous, it cannot
attain both values on the unit sphere, and hence it is constant thereon. \bx

The orientation of $S$ is merely a technicality, the role of which only is to give an
unambigous definition of the map $m$. In our applications of Lemma~\ref{m} we implicitly
assume that an orientation of $S$ has been chosen, an use it without reference.

Our approach to the normal form problem is based on a handy description of the group $\G$.
A triple $(u,v,z)\in V^3$ is called a \emph{Cayley triple} in $V$ if $\{u,v,uv,z\}$ is
an orthonormal set. The set of all Cayley triples in $V$ is denoted by $\mathcal{C}$.

Given $(u,v,z)\in\mathcal{C}$, let $U=\sp\{u,v,uv\}$. For any $x,y\in U$, we have
$\<x,yz\>=\<xy,z\>=0$, as well as $\<z,xz\>=0$. Hence the vector space $V$ decomposes into
an orthogonal direct sum 
$V=U\oplus \sp\{z\}\oplus Uz$. Since $\mathrm{R}_z$ is an iso\-met\-ry on $z^\perp$,
$(uz,vz,(uv)z)$ will be an orthonormal basis for $Uz$. To summarise, this means that every
Cayley triple $c=(u,v,z)$ determines an orthonormal basis
$\underline{b}_c=(u,v,uv,z,uz,vz,(uv)z)$ for $V$. Given a linear endomorphism $\d$ of $V$,
we denote by $[\d]_c$ the matrix of $\d$ with respect to $\underline{b}_c$. 

A group action is said to be \emph{simply transitive} if it is transitive and all stabiliser
subgroups are trivial.

\begin{prop}\label{trans} (See also \cite[11.16]{cpp}.)
  The group $\G$ acts simply transitively on $\mathcal{C}$ by
  $g\!\cdot\!(u,v,z)=(g(u),g(v),g(z))$.
\end{prop}

\textit{Proof:} Clearly, the above expression defines a group action. 
If $c\in\mathcal{C}$, $g\in\G$ and $g\cdot c=c$ then $g(b)=b$ for all
$b\in\underline{b}_c$, and hence $g=\I_V$. So the
stabiliser of $c\in\mathcal{C}$ is trivial.

For transitivity, one must show that the bases given by any two Cayley triples have the
same multiplication table. Note that any permutation of a Cayley triple
$c=(u,v,z)\in\mathcal{C}$ is again a Cayley triple, and that 
$(x,y,z),(xz,y,z)\in\mathcal{C}$ for all ortho\-normal pairs $x,y\in U$.
Therefore, $(xz)z=-x$ and $(xz)(yz)=\(y(xz)\)z=\((xy)z\)z=-xy$. Using this, and Lemma
\ref{vm}, one readily verifies that the multiplication in $V$ in the basis
$\underline{b}_c$ is given by Table~\ref{Vmult}.

\begin{table}[htb] 
  \centering
  \begin{tabular}{r|ccc|c|ccc}
    $\cdot$&$u$&$v$&$uv$       & $z$    &$uz$&$vz$&$(uv)z$\\
    \hline
    $u$ & $0$&$uv$&$-v$        & $uz$   & $-z$&$-(uv)z$&$vz$\\
    $v$ & $-uv$&$0$&$u$        & $vz$   & $(uv)z$&$-z$&$-uz$\\
    $uv$& $v$&$-u$ &$0$        & $(uv)z$& $-vz$&$uz$&$-z$\\
    \hline
    $z$ & $-uz$&$-vz$&$-(uv)z$ & $0$    & $u$&$v$&$uv$\\
    \hline
    $uz$& $z$&$-(uv)z$&$vz$    & $-u$   & $0$&$-uv$&$v$\\
    $vz$& $(uv)z$&$z$&$-uz$    & $-v$   & $uv$&$0$&$-u$\\
    $(uv)z$&$-vz$&$uz$&$z$     & $-uv$  & $-v$&$u$&$0$
  \end{tabular}
  \caption{Multiplication in $(V,\pi)$} \label{Vmult}
\end{table}

The structure constants of $(V,\pi)$ with respect to $\underline{b}_c$ are
independent of the choice of $c\in\mathcal{C}$, and hence for any $c,c'\in\mathcal{C}$
there exists an automorphism $g\in\G$ such that $g\cdot c=c'$. 
Therefore the group action is transitive. \bx

Note that Proposition \ref{trans} implies that $\G$ acts transitively on the
set of orthonormal pairs in $V$, and that the stabiliser of an orthonormal pair $(u,v)$
acts simply transitively on the unit sphere in $\{u,v,uv\}^\perp\subset V$.

Fixing some Cayley triple $s\in\mathcal{C}$, we obtain a bijection
$\mathfrak{t}:\G\rightarrow\mathcal{C},\; g\mapsto g\!\cdot\!s$.
If $g_c=\mathfrak{t}^{-1}(c)$, then $[g_c^{-1}\d g_c]_s=[\d]_{g_c\cdot s}=[\d]_c$.
The normal form problem for (\ref{nform}) can now be rephrased as to find a map
$N:\mathrm{Sym}(V)\rightarrow \mathrm{Sym}(V)$ with the following properties:
\begin{itemize}
\item[\textit{(i)}] $[N(\d)]_s=[\d]_c$ for some $c\in\mathcal{C}$.
\item[\textit{(ii)}] $N(\d)=N(\d')$ whenever there exists a $g\in\G$ such that $\d=g^{-1}{\d'}g$.
\end{itemize}
Then $N(\mathrm{Sym}(V))$ will be a cross-section for $\mathrm{Sym}(V)\slash\G$, and $N(\d)$ the
normal form of $\d$. In other words, for every $\d\in \mathrm{Sym}(V)$, one wants to construct a
non-empty set $\Gamma(\d)\subset\mathcal{C}$ of Cayley triples such that
$[\d]_c=[\d]_{c'}$ for all $c,c'\in\Gamma(\d)$. The set $\Gamma(\d)$ must be chosen only
using properties of $\d$ which are invariant under conjugation with $\G$. A normal form
map then is defined by $[N(\d)]_s=[\d]_c$ where $c\in\Gamma(\d)$.
This construction is carried out (in a somewhat informal way) in
the next section, Section~\ref{normal}.

As $\G\subset \mathrm{O}(V)$, all properties of $\d$ as a linear operator on a
Euclidean space is preserved under conjugation with elements in $\G$. In
particular, the set $\pd=\{(\l,\dim\E_\l(\d))\mid \l \mbox{ is an eigenvalue of } \d \}$ of
eigenpairs of $\d$ is an invariant for its orbit under $\G$. Hence the normal form problem
may be  solved separately for each possible set of eigenpairs.
We distinguish 15 essentially distinct types of sets, determined by the number of
eigenspaces, and their dimensions: 
\begin{itemize}
\item[1:] (7).
\item[2:] (1,6), (2,5), (3,4).
\item[3:] (1,1,5), (1,2,4), (1,3,3), (2,2,3).
\item[4:] (1,1,1,4), (1,1,2,3), (1,2,2,2).
\item[5:] (1,1,1,1,3), (1,1,1,2,2).
\item[6:] (1,1,1,1,1,2).
\item[7:] (1,1,1,1,1,1,1).
\end{itemize}
The different types are treated separately in Section~\ref{normal}.

\section{Normal forms}\label{normal}

The map $\mathcal{s}:\mathrm{Sym}(V)\to \mathrm{Sym}(\R^7),\:\d\mapsto[\d]_s$ is a bijection. The
$\G$-action on $\mathrm{Sym}(V)$ defines a $\G$-action on $\mathrm{Sym}(\R^7)$
by $\mathcal{s}(\d)\cdot g=\mathcal{s}(\d\cdot g)$. A subset $\mathcal{N}\subset
\mathrm{Sym}(\R^7)$ is a cross-section for $\mathrm{Sym}(\R^7)\slash\G$ if and only if
$\mathcal{s}^{-1}(\mathcal{N})$ is a cross-section for $\mathrm{Sym}(V)\slash\G$. For sets $\p$ of
suitable form we define $\mathrm{Sym}_\p=\mathrm{Sym}_\p(\R^7)=\{\d\in \mathrm{Sym}(\R^7)
\mid \pd=\p\}$. In Propositions~\ref{p7}--\ref{p111duo}, cross-sections for
$\mathrm{Sym}_\p\slash\G$ are given for all possible $\p$. The preimage of their union
under $\mathcal{s}$ gives the desired cross-section for $\mathrm{Sym}(V)\slash\G$.

Given $\a,\b,\gamma\in\R$ we write 
\begin{align*}
  R_\a&=\begin{pmatrix}\cos\a&-\sin\a\\ \sin\a&\cos\a\end{pmatrix},\\[0.1cm]
  S_{\a,\b}&=
  \begin{pmatrix}
    \begin{bmatrix}\cos\a\\\sin\a\end{bmatrix} &
      \cos\b\begin{bmatrix}-\sin\a\\\cos\a\end{bmatrix} &
      -\sin\b\begin{bmatrix}-\sin\a\\\cos\a\end{bmatrix} \\
      0 & \sin\b & \cos\b
  \end{pmatrix} = \\[0.1cm]
  &=
  \begin{pmatrix}
    \cos\a & -\cos\b\sin\a & \sin\b\sin\a \\
    \sin\a & \cos\b\cos\a & -\sin\b\cos\a \\
    0 & \sin\b & \cos\b
  \end{pmatrix} \mbox{ and}\\[0.1cm]
  T_{\a,\b,\gamma}&=
  \begin{pmatrix}
    \cos\a &0& -\sin\a\cos\b & \sin\a\sin\b\cos\gamma  & -\sin\a\sin\b\sin\gamma \\
    0      &1& 0&0&0 \\
    \sin\a &0& \cos\a\cos\b  & -\cos\a\sin\b\cos\gamma & \cos\a\sin\b\sin\gamma \\
    0      &0& \sin\b        & \cos\b\cos\gamma        & -\cos\b\sin\gamma \\
    0      &0& 0             & \sin\gamma              & \cos\gamma
  \end{pmatrix}.
\end{align*}
These matrices are used to describe the normal forms, and have certain
geo\-met\-ric interpretations. The matrix $R_\a$ is simply the matrix of rotation in $\R^2$
with the angle $\a$. The matrices of the form $S_{\a,\b}$ constitute the
subset of $\mathrm{SO}_3(\R)$ defined by the property that the image of the first standard
basis vector $e_1$ is orthogonal to $e_3$. Its significance is mainly due to the fact that
the set
\begin{equation*}
  \left\{S_{\a,\b}^{-1}\begin{pmatrix}\l\\&\mu\I_2\end{pmatrix}S_{\a,\b}
    \Big| (\a,\b)\in\left(]0,\pi[\times[0,\pi[\right)\cup\{(0,0)\}\right\}
\end{equation*}
parametrises $\mathrm{Sym}_{\{(\l,1),(\mu,2)\}}$. As for the matrices $T_{\a,\b,\gamma}$,
they make up the set of $\mathrm{SO}_3(\R)$-matrices fixing $e_2$ and mapping $e_1$ into
the orthogonal complement of $e_4$ and $e_5$, and $e_3$ into the orthogonal complement of
$e_5$.

For simplicity, we denote the endomorphism 
\begin{equation*}
  R_\a\otimes\I_V=
  \begin{pmatrix}
    \cos\a\I_V & -\sin\a\I_V \\
    \sin\a\I_V & \cos\a\I_V
  \end{pmatrix}
\end{equation*}
of $V^2$ briefly by $R_\a$. In this notation,
\begin{equation*}
  R_\a\begin{pmatrix}e\\f\end{pmatrix}=
  \begin{pmatrix}\cos\a e-\sin\a f \\ 
    \sin\a e + \cos\a f\end{pmatrix}\mbox{ for } e,f\in V.
\end{equation*}
Throughout this section, $c$ denotes a Cayley triple $(u,v,z)\in\mathcal{C}$.

When solving the normal form problem for (\ref{nform}), we often need to consider the set
$\mathrm{Sym}_\p(\R^k)=\{ \d\in \mathrm{Sym}(\R^k) \mid \pd=\p \}$ of symmetric $k\!\times\!k$-matrices
sharing the same set $\p$ of eigenpairs. In Section~\ref{parametrisering}, $\mathrm{Sym}_\p(\R^k)$
is parametrised for arbitrary $k$ and $\p$. This parametrisation is used without reference
throughout the present section. In most cases, the parametrisation is fairly
uncomplicated, and is then written out explicitly. Otherwise, it is given in terms of the
function $\mathfrak{c}\circ\mathfrak{r}:\mathfrak{K}(D)\rightarrow \mathrm{Sym}_\p$ described in
Section~\ref{parametrisering}, where $D$ is a diagonal matrix such that
$\p(D)=\p$. We use the notation
$\mathfrak{K}'(D)=\left\{\tau\in\mathfrak{K}(D) \mid \tau^{1}_{1}\neq\frac{\pi}{2}\right\}$.

\subsection{Types (7), (1,6), (2,5) and (1,1,5)}

These four types are trivial. First, if $\d$ is of type (7), then there is only one eigenpair,
$(\l,7)$. Any choice of a Cayley triple $c$ will give the matrix $[\d]_c=\l\I_7$. 

If $\d$ is of type (1,6), then we may choose any $c=(u,v,z)\in\mathcal{C}$ such that $u$
belongs to the eigenspace of dimension 1. For the type (2,5), any orthonormal basis
$(u,v)$ for the two-dimensional eigenspace can be extended to a Cayley triple
$c=(u,v,z)$. Finally, if $\pd=\{(\l,1),(\mu,1),(\nu,5)\}$ is the set of eigenpairs of $\d$
(this is the type (1,1,5)), then $c\in\mathcal{C}$ may be chosen such that $u\in\E_\l(\d)$
and $v\in\E_\mu(\d)$, and hence $z\in\E_\nu(\d)$.

The matrices obtained will be
\begin{alignat*}{2}
  &[\d]_c=
  \begin{pmatrix}
    \l\\
    &\mu\I_6
  \end{pmatrix}
  &&\mbox{for (1,6),} \\
  &[\d]_c=
  \begin{pmatrix}
    \l\I_2\\
    &\mu\I_5
  \end{pmatrix}
  &&\mbox{for (2,5) and} \\
  &[\d]_c=
  \begin{pmatrix}
    \l\\
    &\mu\\
    &&\nu\I_5
  \end{pmatrix}
  &&\mbox{for (1,1,5)}
\end{alignat*}
where certainly all the parameters $\l$, $\mu$ and $\nu$ are determined by the map $\d$
itself, independent of the basis $\underline{b}_c$.
Thus we have proved the following proposition:

\begin{prop}\label{p7}
  \begin{enumerate}
  \item The set $\{\l\I_7\}$ is a cross-section for $\mathrm{Sym}_{\{(\l,7)\}}\slash\G$.
  \item The set $\left\{\left(\begin{smallmatrix}\l\\&\mu\I_6\end{smallmatrix}\right)\right\}$ is a
    cross-section for $\mathrm{Sym}_{\{(\l,1),(\mu,6)\}}\slash\G$. 
  \item The set
    $\left\{\left(\!\begin{smallmatrix}\l\I_2\\&\mu\I_5\end{smallmatrix}\!\right)\right\}$
    is a cross-section for $\mathrm{Sym}_{\{(\l,2),(\mu,5)\}}\slash\G$.
  \item The set
  $\left\{\left(\!\begin{smallmatrix}\l\\&\mu\\&&\nu\I_5\end{smallmatrix}\!\right)\right\}$
  is a cross-section for $\mathrm{Sym}_{\{(\l,1),(\mu,1),(\nu,5)\}}\slash\G$.
  \end{enumerate}
\end{prop}

\subsection{Types (1,1,1,4), (1,2,4) and (3,4)}

The endomorphisms $\d\in \mathrm{Sym}(V)$ of these types have the common property of having a
4-dimensional eigenspace, which we denote by $U$. Suppose that $u,v\in U^\perp$ are
ortho\-normal. Equipping $S=U^\perp$ with some orientation,
we have $uv=m(x)$ for some $x\in S$ of unit length. From Lemma~\ref{m} follows that the
number $\<S,uv\>$ is independent of the choices of $u$ and $v$. Taking
$z\in\{u,v,uv\}^\perp$ such that $U^\perp\subset\sp\{u,v,uv,z\}$ gives a Cayley triple
$c=(u,v,z)\in\mathcal{C}$ for which $uz,vz,(uv)z\in U$ and
$\left(\begin{smallmatrix}uv\\z\end{smallmatrix}\right)=
R^{-1}_\a\left(\begin{smallmatrix}e\\f\end{smallmatrix}\right)$
where $e\in U^\perp$, $f\in U$ and $\a\in[0,\pi[$. By choosing the sign of $z$, it is
possible to obtain $\a\in\left[0,\frac{\pi}{2}\right]$ (replacing $z$ with $-z$
changes $\a$ to $\pi-\a$). Since now $\cos\a=\<S,uv\>$, the angle $\a$ does not depend on the
choices we made in the construction of $c=(u,v,z)$. 

Note that the unit vectors $u$ and $v$ are chosen freely in $U^\perp$. For simplicity, we
want them to be eigenvectors, and in the case (1,2,4) to lie in the eigenspace of
dimension 2. Our construction yields the following results for the different types of
endomorphisms.

\begin{prop}\label{p1114}
  If $\p=\{(\l_1,1),(\l_2,1),(\l_3,1),(\mu,4)\}$, then the set 
  \begin{align*}
    \left\{
    \begin{pmatrix}
      \l_1\\
      &\l_2\\
      &&D\\
      &&&\mu\I_3
    \end{pmatrix}\bigg|
    D=R_\a^{-1}\begin{pmatrix}\l_3\\&\mu\end{pmatrix}R_\a,\;
      \a\in\left[0,\frac{\pi}{2}\right]\right\}
  \end{align*}
  is a cross-section for $\mathrm{Sym}_\p\slash\G$.  
\end{prop}

\begin{prop}\label{p124}
  If $\p=\{(\l,1),(\mu,2),(\nu,4)\}$, then the set
  \begin{align*}
    \left\{
    \begin{pmatrix}
      \mu\I_2\\
      &D\\
      &&\nu\I_3
    \end{pmatrix}\bigg|
    D=R_\a^{-1}\begin{pmatrix}\l\\&\nu\end{pmatrix}R_\a,\;
      \a\in\left[0,\frac{\pi}{2}\right]\right\}
  \end{align*}
  is a cross-section for $\mathrm{Sym}_\p\slash\G$.  
\end{prop}

\begin{prop}\label{p34}
  If $\p=\{(\l,3),(\mu,4)\}$, then the set
  \begin{align*}
    \left\{
    \begin{pmatrix}
      \l\I_2\\
      &D\\
      &&\mu\I_3
    \end{pmatrix}\bigg|
    D=R_\a^{-1}\begin{pmatrix}\l\\&\mu\end{pmatrix}R_\a,\;
      \a\in\left[0,\frac{\pi}{2}\right]\right\}
  \end{align*}
  is a cross-section for $\mathrm{Sym}_\p\slash\G$.  
\end{prop}

\subsection{Types (1,3,3), (2,2,3) and (1,1,2,3)}

Although these classes of endomorphisms differ from each other in im\-port\-ant aspects,
the strategy for dealing with them is roughly the same. Therefore they are treated together.

\begin{lma}\label{lbig3}
  Let $V=X \oplus Y\oplus Z$ be a decomposition of $V$ into pairwise ortho\-gonal non-trivial
  subspaces, where $\dim X=3$ and $\dim Y\geqslant2$. There exist $c\in\mathcal{C}$ and
  $\a\in\R$ such that $u,v\in X$,
  $\left(\begin{smallmatrix}uv\\z\end{smallmatrix}\right)=
    R^{-1}_\a\left(\begin{smallmatrix}e\\f\end{smallmatrix}\right)$
  for some unit vectors $(e,f)\in X\times Y$, and
\begin{equation*}
  uz\in
  \begin{cases}
    Y&\mbox{if }\: (\dim Y,\dim Z)=(3,1),\\
    Z&\mbox{if }\: (\dim Y,\dim Z)=(2,2).
  \end{cases}
\end{equation*}
\end{lma}

\textit{Proof:} If $X$ is closed under $\pi$, let $z$ be any unit vector in $Y$. Now
$\mathrm{R}_z$ maps $X$ bijectively onto $(Y\cap z^\perp)\oplus Z$. In the case
$(\dim Y,\dim Z)=(3,1)$ we choose $u,v\in X$ such that
$u\in \mathrm{R}^{-1}_z(Y\cap z^\perp)$, and in the case $(\dim Y,\dim Z)=(2,2)$ such that
$u\in \mathrm{R}^{-1}_z(Z)$. 

Now suppose $X$ is not closed under $\pi$. Set $S=X$. Non-closedness of $S$ implies that
the linear map $P_{Y\oplus Z}m:S\to Y\!\oplus\! Z$ is injective. Thus we have
$\dim(P_{Y\oplus Z}m(S))=3$, and hence $Y\cap (P_{Y\oplus Z}m)(S)$ is non-empty. Let
$x\in(P_{Y\oplus Z}m)^{-1}(Y)$ be a unit vector. We get $m(x)=\cos\a e-\sin\a f$ for some $\a\in\R$,
$(e,f)\in X\times Y$. Set $z=\sin\a e+\cos\a f$. 

The subspace $\mathrm{R}_z(X\cap x^\perp)\subset(Y\cap m(x)^\perp)\oplus Z$ is 2-dimensional.
If $\dim Y=3$ then $\dim(Y\cap m(x)^\perp)=2$, so
$\mathrm{R}_z(X\cap x^\perp)\cap(Y\cap m(x)^\perp)\neq0$. Hence there exists a unit vector
$u\in X\cap x^\perp$ such that $uz\in Y$. If $\dim Z=2$, by the same argument there exists
a unit vector $u\in X\cap x^\perp$ such that $uz\in Z$. In both cases, we may add $v\in X$ such that
$(u,v,x)$ becomes a positively oriented orthonormal basis for $S$. Since now $uv=m(x)$ and
$\<m(x),z\>=0$, we have a Cayley triple $(u,v,z)\in\mathcal{C}$ with the desired
properties. \bx

\begin{prop}\label{p133}
  Let $\p=\{(\l,1),(\mu,3),(\nu,3)\}$. A cross-section for $\mathrm{Sym}_\p\slash\G$ is given by
  \begin{equation*}
    \left\{
     \begin{pmatrix}
       \mu\I_2\\
       &A\\
       &&\nu\\
       &&&B
     \end{pmatrix} \left|
     \begin{array}{l}
       A=R_\a^{-1} \begin{pmatrix}\mu\\&\nu\end{pmatrix} R_{\a},\;
	 B=R_\theta^{-1} \begin{pmatrix}\nu\\&\l\end{pmatrix} R_\theta\\
	   \rule{0cm}{5mm}
	   (\a,\theta)\in\left(\left]0,\frac{\pi}{2}\right] \!\times\!
	 \left[0,\frac{\pi}{2}\right]\right) 
	 \cup\{(0,0)\}
     \end{array}\right.\right\}.
  \end{equation*}
\end{prop}

\textit{Proof:} From Lemma~\ref{lbig3} follows (set $X=\E_\mu$, $Y=\E_\nu$ and $Z=\E_\l$)
that for each $\d\in \mathrm{Sym}_\p$ there exists a Cayley triple
$c'=(u',v',z')\in\mathcal{C}$ such that
\begin{equation}\label{m133}
[\d]_{c'}=
  \begin{pmatrix}
       \mu\I_2\\
       &A\\
       &&\nu\\
       &&&B
     \end{pmatrix}
\end{equation}
with $A$ and $B$ as in the proposition, and $\a,\theta\in[0,\pi[$. If $\E_\mu$
is closed under $\pi$ (that is if $\a=0$) the restricted and co-restricted map
$\mathrm{R}_{z'}:\E_\mu\to(\E_\nu\cap z'^\perp)\oplus\E_\l$ is a bijection. Thus there exist
orthonormal vectors $u,v\in\E_\mu$ such that $uz',vz'\in\E_\nu$, and consequently
$(uv)z'\in\E_\l$. On setting $z=z'$, the triple $c=(u,v,z)$ becomes a Cayley triple, for
which $[\d]_c$ has the form stated in the proposition, with $\a=\theta=0$.

Assume $\E_\mu$ is not closed under $\pi$. Set $z=z'$ if $\a\leqslant\frac{\pi}{2}$ in
(\ref{m133}) and $z=-z'$ if $\a>\frac{\pi}{2}$.
If $\theta>\frac{\pi}{2}$, then let $(u,v)=(-u',-v')$, otherwise take
$(u,v)=(u',v')$. This gives us a Cayley triple $c=(u,v,z)$ for which $[\d]_c$ has the form
given in the proposition, with $\a\neq0$.
If we set $S=\E_\mu$, we will have
$\cos\a=\<\E_\mu,m(x)\>$ and $(uv)z=\frac{1}{\<\E_\mu,m(x)\>}xm(x)$, where
$x\in\E_\mu\cap\{u,v\}^\perp$ is the unique vector such that $uv=m(x)$. From Lemma~\ref{m}
follows that $\cos\a$ and $(uv)z$ are the same for all possible choices of
$c\in\mathcal{C}$ for which $[\d]_c$ has the above form. Therefore different
pairs $(\a,\theta)\in]0,\frac{\pi}{2}]\times[0,\frac{\pi}{2}]$ cannot correspond to
endomorphisms within the same orbit of the group action (\ref{nform}). \bx

The case (2,2,3) is analogous to the case (1,3,3). By the same technique as for
Proposition~\ref{p133}, using Lemma~\ref{lbig3} with $(X,Y,Z)=(\E_\nu,\E_\l,\E_\mu)$, the
following proposition is proved.

\begin{prop}\label{p223}
  If $\p=\{(\l,2),(\mu,2),(\nu,3)\}$, then the set
  \begin{equation*}
    \left\{
    \begin{pmatrix}
      \nu\I_2\\
      &A\\
      &&\mu\\
      &&&B
    \end{pmatrix} \left|
    \begin{array}{l}
      A=R_\a^{-1} \begin{pmatrix}\nu\\&\l\end{pmatrix} R_{\a},\;
	B=R_\theta^{-1} \begin{pmatrix}\l\\&\mu\end{pmatrix} R_{\theta}\\
	  \rule{0cm}{5mm}
	  (\a,\theta)\in\left(\left]0,\frac{\pi}{2}\right] \!\times\!
	\left[0,\frac{\pi}{2}\right] \right) \cup \{(0,0)\}
    \end{array}
    \right.\right\}
  \end{equation*}
  is a cross-section for $\mathrm{Sym}_\p\slash\G$.
\end{prop}

The case (1,1,2,3) is somewhat different from the two previous cases. The result reads as
follows.

\begin{prop}\label{p1123}
  If $\p=\{(\kappa,1),(\l,1),(\mu,2),(\nu,3)\}$ then
  \begin{equation*}
    \left\{
    \begin{pmatrix}
      \nu\I_2\\ &A \\ &&B
    \end{pmatrix} \left|
    \begin{array}{l}
      A=R_\a^{-1}\left(\!\begin{smallmatrix} \nu \\ &\mu \end{smallmatrix}\!\right)R_\a,\;
      B=S_{\theta,\phi}^{-1}\left(\!\begin{smallmatrix}\mu\\&\l\\&&\kappa\end{smallmatrix}\!\right)S_{\theta,\phi},\\
      \rule{0cm}{1cm}
	\begin{aligned}
	  (\a,\theta,\phi)\in&\left]0,\frac{\pi}{2}\right]\!\times\!
	  \left[0,\frac{\pi}{2}\right]\!\times\!\left]0,\frac{\pi}{2}\right] \\ 
	  \rule{0cm}{5mm}
	  &\cup\left(\left]0,\frac{\pi}{2}\right]\!\times\!\{(0,0)\}\right)\cup\{(0,0,0)\}
	\end{aligned}
    \end{array}
    \right.\right\}
  \end{equation*}
  is a cross-section for $\mathrm{Sym}_\p\slash\G$.
\end{prop}

\textit{Proof:} If $\E_\nu$ is closed under $\pi$ then, analogous to the case (1,3,3), one
shows the existence of a Cayley triple $c=(u,v,z)\in\mathcal{C}$ for which
$u,v,uv\in\E_\nu,\:z,uz\in\E_\mu$, $vz\in\E_\l$ and $(uv)z\in\E_\kappa$. This is the case
$\a=\theta=\phi=0$.

Suppose $\E_\nu$ is not closed under $\pi$. Set $S=\E_\nu$. Now 
$P_{S^\perp}m:S\to S^\perp$ is injective and hence, for dimension reasons, there exists a
unit vector $x\in S=\E_\nu$ such that $m(x)\in(\E_\nu\oplus\E_\mu)\minus\E_\nu$. In this
case $\{x,m(x)\}$ is linearly independent, so we may chose $z\in\sp\{x,m(x)\}$ to be a
unit vector orthogonal to $m(x)$. By Lemma~\ref{m}:\ref{m3}, the vector $x(m(x))$ is
independent of the choice of $x\in\E_\nu$. We get two subcases:

First case: $x(m(x))\in\E_\kappa$. This means that
$z(m(x))=\frac{1}{\<\E_\nu^\perp,m(x)\>}x(m(x))\in\E_\kappa$. It follows that
$\mathrm{R}_z(\E_\nu\cap x^\perp)=(\E_\mu\cap m(x)^\perp)\oplus\E_\l$, whence there exist orthonormal
$u,v\in\E_\nu\cap x^\perp$ such that $uz\in\E_\mu$ and $vz\in\E_\l$. After choosing the
sign of $u$, we arrive at the case
$(\a,\theta,\phi)\in\left]0,\frac{\pi}{2}\right]\!\times\!\{(0,0)\}$ of the proposition.

Second case: $x(m(x))\not\in\E_\kappa$. Here $z(m(x))\not\in\E_\kappa$, and
$\dim(\mathrm{R}_z(\E_\nu\cap x^\perp)\cap((\E_\mu\cap m(x)^\perp)\oplus\E_\l))=1$. Hence there
exists a unit vector $u\in\E_\nu\cap x^\perp$ (unique up to change of sign) such that
$uz\in(\E_\mu\cap m(x)^\perp)\oplus\E_\l$. Adding $v\in\E_\nu\cap\{u,x\}^\perp$ of unit length, and
possibly changing the signs of the vectors $u$, $v$ and $z$, we get $[\d]_c$ as in the
proposition with $(\a,\theta,\phi)\in\left]0,\frac{\pi}{2}\right]\!\times\!
\left[0,\frac{\pi}{2}\right]\!\times\!\left]0,\frac{\pi}{2}\right]$.  

To prove the irredundancy of the parametrisation given in the proposition, we first note
that $\cos\a=\<\E_\nu,m(x)\>$, where $x$ is any unit vector in $\E_\nu$. This means that
the angle $\a\in\left[0,\frac{\pi}{2}\right]$ is uniquely determined by $\d$. 
In case $\E_\nu$ is not closed under $\pi$, we also have
$(uv)z=\pm\frac{1}{\<\E_\nu^\perp,m(x)\>}x(m(x))$, and hence
$\cos\phi=\frac{\<\E_\kappa,x(m(x))}{\<\E_\nu^\perp,m(x)\>}$ which is independent of the
choice of $x$. Therefore, $\phi$ is also independent of the choice of the Cayley triple
$c$. Finally, if $\phi\neq0$ then $\sin\theta=\frac{\<\E_\mu,(uv)z\>}{\sin\phi}$, which
implies independence of choice for the angle $\theta$. Hence all the angles are determined by $\d$
itself, and the parametrisation is irredundant. \bx

\subsection{Type (1,2,2,2)}

This is perhaps the most difficult of all 15 cases. We begin by spelling out a rather
obvious fact.

\begin{lma}\label{mm}
  Let $U\subset V$ be a 2-dimensional subspace, and $x\in U^\perp$. Then the function
  $v\mapsto\<U,xv\>$ is constant on the unit sphere in $U$.
\end{lma}

\textit{Proof:} This is a consequence of Lemma~\ref{m}:\ref{m2}. Set
$S=U\oplus\sp\{x\}$. Then $xv=\|x\|m(e)$ for some $e\in S$. Since $\<xv,x\>=0$ we have
$\<U,xv\>=\<S,xv\>=\|x\|\<S,m(e)\>$. Now $e\mapsto\<S,m(e)\>$ is constant on the unit
sphere in $S$, whence $v\mapsto\<U,xv\>$ is constant on the unit sphere in $U$. \bx

Suppose $\d\in \mathrm{Sym}(V),\:\pd=\{(\kappa,1),(\l,2),(\mu,2),(\nu,2)\}$. 
If $S=\E_\kappa\oplus\E_\l$ is closed under $\pi$, take $z\in\E_\mu$ to be a unit
vector. Then $\mathrm{R}_z(\E_\l)\subset(\E_\mu\cap z^\perp)\oplus\E_\nu$, which implies
existence of a unit vector $u\in\E_\l$ for which $uz\in\E_\nu$. Choosing
$v\in\E_\l\cap u^\perp$ appropriately, we get a Cayley triple $c=(u,v,z)$ such that
\begin{align*}
  [\d]_c&=
  \begin{pmatrix}
    \l\I_2\\
    &\kappa\\
    &&\mu\\
    &&&\nu\\
    &&&&B'\\
  \end{pmatrix},\\
  B'&=R_\phi^{-1} \begin{pmatrix}\mu\\&\nu\end{pmatrix}R_\phi,
    \mbox{ with } \phi\in\left[0,\frac{\pi}{2}\right].
\end{align*}
By Lemma~\ref{mm}, the value of $\sin\phi=\<\E_\mu,(uv)z\>$ does not depend on the choice of
$z$. Therefore, the matrix $B'$ is uniquely determined by $\d$.

Suppose $S$ is not closed under $\pi$. Set
$r=P_{S^\perp}m:S\rightarrow\E_\mu\oplus\E_\nu$. By Lemma~\ref{m}:\ref{m2}, we have 
$\dim r(S)=3$, which implies $r(S)\cap\E_\mu\neq\{0\}$. 
We dis\-tin\-guish three different subcases:

First case: $\dim r^{-1}(\E_\mu)=1,\;r(\E_\kappa)\subset\E_\mu$. This implies
$m(\E_\kappa)\subset\E_\kappa\oplus\E_\mu$. Because $S$ is not closed under $\pi$, we have
$m(\E_\kappa)\neq\E_\kappa$. Therefore, the subspace $(\E_\kappa+m(\E_\kappa))\cap
m(\E_\kappa)^\perp$ of $V$ is non-trivial.
Let $z\in(\E_\kappa+m(\E_\kappa))\cap m(\E_\kappa)^\perp$ be a unit vector. As
$\mathrm{R}_z(\E_\l)\subset(\E_\mu\cap z^\perp)\oplus\E_\nu$, there exists a unit vector
$u\in\E_\l$ such that $uz\in\E_\nu$. Adding $v\in\E_\l\cap u^\perp$ with suitable
orientation, and possibly changing the sign of $z$ (this is to obtain
$\b,\phi\leqslant\frac{\pi}{2}$ below), we get $c=(u,v,z)\in\mathcal{C}$ for which
\begin{align*}
  [\d]_c&=
  \begin{pmatrix}
    \l\I_2\\
    &A'\\
    &&\nu\\
    &&&B'
  \end{pmatrix},\\
  A'&=R_\b^{-1}\begin{pmatrix}\kappa\\&\mu\end{pmatrix}R_\b,\;
    B'=R_\phi^{-1}\begin{pmatrix}\mu\\&\nu\end{pmatrix}R_\phi,\;\;
      (\b,\phi)\in\left]0,\frac{\pi}{2}\right]\!\times\!\left[0,\frac{\pi}{2}\right].
\end{align*}
Note that we have $\b\neq0$, since $S$ is not closed under $\pi$.

Second case: $\dim r^{-1}(\E_\mu)=1,\;r(\E_\kappa)\not\subset\E_\mu$. Let
$e\in r^{-1}(\E_\mu),\|e\|=1$. By assumption, $e$ is unique up to change of sign, and
$e\not\in\E_\kappa$. Hence $\<\E_\l,e\>\neq0$ and thereby $\dim(\E_\l\cap e^\perp)=1$.
Take $u\in\E_\l\cap e^\perp$ and $v\in S\cap\{e,u\}^\perp$ to be unit vectors. Both $u$
and $v$ are chosen in 1-dimensional subspaces, and are therefore uniquely determined by
$\d$ up to change of sign. Since $S$ is not closed under $\pi$, the vectors $e$ and $uv$ are
non-proportional, and we may choose
$z\in\sp\{e,uv\}$ to be a unit vector orthogonal to $uv$.
Then $c=(u,v,z)$ is a Cayley triple, and
\begin{align*}
  [\d]_c&=
  \begin{pmatrix}
    \l\\
    &A\\
    &&B
  \end{pmatrix}, \mbox{ where}\\
  A&=S_{\a,\b}^{-1}\begin{pmatrix}\l\\&\kappa\\&&\mu\end{pmatrix}S_{\a,\b},\;
    \a,\b\in\left]0,\pi\right[ \mbox{ and}\\
	B&=S_{\theta,\phi}^{-1}\begin{pmatrix}\mu\\&\nu\I_2\end{pmatrix}S_{\theta,\phi},\;
	  \theta,\phi\in\left[0,\pi\right[.
\end{align*}

\begin{lma} \label{lvinklar}
  The triple $c\in\mathcal{C}$ can be chosen such that either
  \begin{enumerate}
  \item $\a,\b,\theta\in\left]0,\frac{\pi}{2}\right[,\;\phi\in[0,\pi[$ \mbox{ or}
  \item
    $\a,\b,\theta\in\left]0,\frac{\pi}{2}\right],\; \frac{\pi}{2}\in\{\a,\b,\theta\}, \;
	      \phi\in\left[0,\frac{\pi}{2}\right]$ \mbox{ or}
  \item $\a,\b\in\left]0,\frac{\pi}{2}\right],\;\theta=\phi=0$.
  \end{enumerate}
  In this presentation, the tuple $(\a,\b,\theta,\phi)$ is uniquely determined by $\d$.
\end{lma}
We denote by $\mathcal{L}_2$ the set of angles $(\a,\b,\theta,\phi)\in\R^4$ described in
Lemma~\ref{lvinklar}.

\textit{Proof:} In our construction, the vectors $u$, $v$ and $z$ (an thereby all elements
of the basis $\underline{b}_c$) are uniquely determined up to change of sign.
In any case, by choosing the signs of $u$, $v$ and $z$ we can assure
that $\a,\b,\theta\leqslant\frac{\pi}{2}$, but in general not
$\phi\leqslant\frac{\pi}{2}$. If $\{v,uv,uz\}$ contains an eigenvector of $\d$, indeed it
is possible also to get $\phi\leqslant\frac{\pi}{2}$. 
In the situation when $uz\in\E_\mu$ (\ie{} when $\theta=0$) we will have
$vz,(uv)z\in\E_\nu$, and may set $\phi=0$. These are in turn the cases 1, 2 and 3 in the
lemma. \bx

Third case: $\dim(r^{-1}(\E_\mu))=2$.
Since $\E_\l\subset S$ has codimension $1$, it follows that
$r^{-1}(\E_\mu)\cap\E_\l\neq\{0\}$. Hence there is a unit vector $e\in\E_\l$ such that
$r(e)\in\E_\mu$. This implies the existence of unit vectors $u\in\E_\l$, $v\in\E_\kappa$
with the property that $uv=m(e)\in\E_\l\oplus\E_\mu$. Taking $z\in\sp\{e,uv\}\cap uv^\perp$ we get
a Cayley triple $c=(u,v,z)$.

\begin{lma}\label{m1222}
  Given $c\in\mathcal{C}$ as above, $\dim r^{-1}(\E_\mu)=2$ if and only if
  $(uv)z\in\E_\nu$.
\end{lma}

\textit{Proof:} For any unit vector $x\in S$ we have
$(uv)z=\pm\frac{1}{\<S^\perp,m(x)\>}xm(x)$.
This implies, by Lemma~\ref{m}, that $(uv)z\in(S+m(S))^\perp=(S+r(S))^\perp$.
Since $r:S\rightarrow\E_\mu\oplus\E_\nu$ is injective, for dimension reasons we have
$\E_\mu\oplus\E_\nu=r(S)\oplus\sp\{(uv)z)\}$. Certainly, 
$\dim r^{-1}(\E_\mu)=\dim(r(S)\cap\E_\mu)=2$ if and only if $\<\E_\mu,(uv)z\>=0$, that is if
and only if $(uv)z\in\E_\nu$. \bx

We remark that $r(\E_\l)=\E_\mu$ precisely when $\E_\nu=\sp\{vz,(uv)z\}$. 

Let 
\begin{align*}
\mathcal{L}_1&=\{0\}\times\left[0,\frac{\pi}{2}\right]\times
\left\{\frac{\pi}{2}\right\}\times\left[0,\frac{\pi}{2}\right] \\
\intertext{and}
\mathcal{L}_3&=\left\{\frac{\pi}{2}\right\}\times\left]0,\frac{\pi}{2}\right]
\times \left[0,\frac{\pi}{2}\right]\times\{0\}.
\end{align*}
Summarising all cases above, we get the following proposition.

\begin{prop}\label{p1222}
  Let $\p=\{(\kappa,1),(\l,2),(\mu,2),(\nu,2)\}$. A cross-section for $\mathrm{Sym}_\p\slash\G$ is
  given by
  \begin{equation*}
    \left\{
    \begin{pmatrix}
      \l\\
      &A\\
      &&B
    \end{pmatrix} 
    \left|
    \begin{array}{l}
      A=S_{\a,\b}^{-1}\left(\!\begin{smallmatrix}\l\\&\kappa\\&&\mu\end{smallmatrix}\!\right)
      S_{\a,\b}\,,\;
      B =S_{\theta,\phi}^{-1}\left(\!\begin{smallmatrix}\mu\\&\nu\I_2\end{smallmatrix}\!\right)
      S_{\theta,\phi} \\
      \rule{0cm}{5mm}
      (\a,\b,\theta,\phi)\in\mathcal{L}_1\cup\mathcal{L}_2\cup\mathcal{L}_3
    \end{array}
    \right. \right\}
  \end{equation*}
\end{prop}
The set of angles $\mathcal{L}_1$ here covers the case when $S$ is closed under $\pi$
(this is when $\b=0$), and what is referred to as first case above. The sets
$\mathcal{L}_2$ and $\mathcal{L}_3$ correspond to the second and third cases respectively.

\subsection{Type (1,1,1,2,2)}

\begin{prop}\label{p11122}
  Let $\p=\{(\l_1,1),(\l_2,1),(\l_3,1),(\mu,2),(\nu,2)\}$ and let
  $D=\left(\begin{smallmatrix}\mu\rule{1mm}{0mm}\\&\l_3\\&&\nu\I_2\end{smallmatrix}\right)$. Then
  \begin{align*}
    &\left\{
    \begin{pmatrix}
      \l_1\\
      &\l_2\\
      &&\l_3\\
      &&&\mu\\
      &&&&B
    \end{pmatrix}
    \left|
    \begin{array}{l}
      B=S^{-1}_{\a,\b}\begin{pmatrix}\mu\\&\nu\I_2\end{pmatrix}S_{\a,\b},\\
      \rule{0cm}{5mm}
      (\a,\b)\in\left(\left]0,\frac{\pi}{2}\right]\times\left[0,\frac{\pi}{2}\right]\right)
      \cup\{(0,0)\}
    \end{array}
    \right.\right\} \dot\cup \\[0.1cm]
    &\left\{
    \left(\!
    \begin{array}{cccc|c}
      \l_1&&&&0\\
      &\l_2&&&0\\
      && a&&\mathfrak{a}^\ast\\
      && &\mu&0\\ 
      \hline
      0&0& \mathfrak{a}&0&A
    \end{array}
    \!\right)
    \left|
    \begin{array}{l}
      a\in\R\!\minus\!\{0\},\:\mathfrak{a}\in\R^3,\:A\!\in\!\mathrm{Sym}(\R^3),\\ 
      \rule{0cm}{6mm}
      \begin{pmatrix}
        a&\mathfrak{a}^\ast\\
	\mathfrak{a}&A
      \end{pmatrix}\in(\mathfrak{c}\circ\mathfrak{r})(\mathfrak{K}'_{1,2,3}(D))
    \end{array}
    \right.\right\} \dot\cup \\[0.1cm]
    &\left\{
    \begin{pmatrix}
      \l_1\\
      &\l_2\\
      &&C
    \end{pmatrix}
    \left|
    \begin{array}{l}
      C=T^{-1}_{\a,\b,\gamma}
      \left(\!\begin{smallmatrix}
	\nu\I_2\\
	&\l_3\\
	&&\mu\I_2
      \end{smallmatrix}\!\right) T_{\a,\b,\gamma},\\
      \rule{0cm}{5mm}
      (\a,\b,\gamma)\in\left[0,\frac{\pi}{2}\right[
	  \!\times\!\left[0,\frac{\pi}{2}\right]\!\times\!\left[0,\frac{\pi}{2}\right]
    \end{array}
    \right.\right\}
  \end{align*}
  is a cross-section for $\mathrm{Sym}_\p\slash\G$.
\end{prop}

\textit{Proof:} Suppose $\d\in \mathrm{Sym}_\p$. Take $(u,v)\in\E_{\l_1}\!\times\!\E_{\l_2}$. We
distinguish three subcases:

First case: $uv\in\E_{\l_3}$. Here let $z\in\E_\mu$. With a suitable choice of signs for
$u$, $v$ and $z$, the matrix $[\d]_c$ will belong to the first set given in the proposition.
From Lemma~\ref{mm} follows that $\cos\a=\<\E_\mu,uz\>$ and $-\cos\b\sin\a=\<\E_\mu,vz\>$
do not depend on the choice of $z\in\E_\mu$. Since alterations of the vectors $u$ and $v$
can only change $\a$ and  $\b$ to $\pi-\a$ and  $\pi-\b$ respectively (and thereby give
rise to matrices that are either unchanged or do not belong to the set), the given matrices belong
to different orbits of (\ref{nform}). 

Second case: $\<\E_\mu,uv\>\neq0$. Taking $z\in\E_\mu\cap uv^\perp$, the Cayley triple
$c\in\mathcal{C}$ is uniquely determined up to change of signs. This gives the second set
in the proposition.

Third case: $\<\E_\mu,uv\>=0,\; uv\not\in\E_{\l_3}$. Since in this case
$\<\E_\nu,uv\>\neq0$, a unit vector $z\in\E_\nu\cap uv^\perp$ can be chosen in only two
ways. The matrices $[\d]_c$ obtained here constitute the third set. \bx

\subsection{Types (1,1,1,1,1,1,1), (1,1,1,1,1,2) and (1,1,1,1,3)}

These three sorts of matrices can be treated simultaneously. Let
\begin{equation*}
  \a\in\R,\;\;
  X=
  \begin{pmatrix}
    x_1&x_2&x_3&x_4\\
    X_1&X_2&X_3&X_4
  \end{pmatrix}\in\R^{4\times4},\;\; x_i\in\R,\;X_i\in\R^3,
\end{equation*}
and define
\begin{multline*}
  U_{X,\a}\!=\!
  \left(\!\!
  \begin{array}{ccccc}
    x_1\!\begin{bmatrix}-\sin\a\\\cos\a\end{bmatrix} &
    \begin{bmatrix}\cos\a\\\sin\a\end{bmatrix} &
      x_2\!\begin{bmatrix}-\sin\a\\\cos\a\end{bmatrix} &
      x_3\!\begin{bmatrix}-\sin\a\\\cos\a\end{bmatrix} &
      x_4\!\begin{bmatrix}-\sin\a\\\cos\a\end{bmatrix}  \\
      X_1 & 0 & X_2 & X_3 & X_4
  \end{array}
  \!\!\right) \\[0.1cm]
  =
  \begin{pmatrix}
    -x_1\sin\a & \cos\a & -x_2\sin\a & -x_3\sin\a & -x_4\sin\a \\
    x_1\cos\a & \sin\a & x_2\cos\a & x_3\cos\a & x_4\cos\a \\
    X_1 & 0 & X_2 & X_3 & X_4    
  \end{pmatrix}
  \in\R^{5\times5}.
\end{multline*}
Given $1\leqslant j \leqslant k \leqslant 6$ and $\l_j,\ldots,\l_k,\mu\in\R$, set
\begin{equation*}
  D^{(j)}_k=
  \begin{pmatrix}
    \l_j\\&\ddots\\&&\l_k\\&&&\mu\I_{7-k}
  \end{pmatrix}.
\end{equation*}

\begin{prop}\label{p111duo}
  Suppose $\p=\{(\l_i,1),(\mu,7-k)\}_{i\in\underline{k}}$,  $k\in\{4,5,6\}$.
  A cross-section for $\mathrm{Sym}_\p\slash\G$ is given by 
  \begin{align*}
    &\left\{
    \left(\!
    \begin{array}{cccc|c}
      \l_1&&&&0\\
      &\l_2&&&0\\
      && a&&\mathfrak{a}^\ast\\
      && &\l_3&0\\
      \hline
      0&0& \mathfrak{a}&0&A
    \end{array}\!\right) \left|
    \begin{array}{l}
      a\in\R,\;\mathfrak{a}\in\R^3,\:A\in \mathrm{Sym}(\R^3) \\
      \rule{0cm}{6mm}
      \begin{pmatrix}
	a&\mathfrak{a}^\ast\\
	\mathfrak{a}&A
      \end{pmatrix}
      \in(\mathfrak{c}\circ\mathfrak{r})(\mathfrak{K}_{1,2,3}(D^{(4)}_k))
    \end{array} \right.\right\} \dot\cup \\[0.1cm]
&\left\{
    \begin{pmatrix}
      \l_1\\
      &\l_2\\
      &&B
    \end{pmatrix} \left|
    \begin{array}{l}
      B=U_{X,\a}^{-1} D^{(3)}_k\, U_{X,\a} \\
      \rule{0cm}{4mm}
      (\a,X)\in
      \left]0,\frac{\pi}{2}\right[\times\mathfrak{r}(\mathfrak{K}'_{1,2}(D^{(3)}_k)) \\
      \rule{0cm}{4mm}
      \hfill\cup\left\{\frac{\pi}{2}\right\}\times\mathfrak{r}(\mathfrak{K}'_{1,2,3}(D^{(3)}_k))
    \end{array}
    \right.\right\}.
  \end{align*}  
\end{prop}

\textit{Proof:}
Taking $(u,v)\in\E_{\l_1}\!\times\!\E_{\l_2}$ we get two subcases.
Firstly, if $\<\E_{\l_3},uv\>=0$ then we let $z\in\E_{\l_3}$. This gives
$c=(u,v,z)\in\mathcal{C}$ for which $[\d]_c$ belongs to the first set in the proposition.

Secondly, if $\<\E_{\l_3},uv\>\neq0$, then we set $E=\E_{\l_3}\oplus\E_{\l_4}$ and choose
$z$ in the 1-dimensional subspace $E\cap uv^\perp$.
Writing $w=(P_E(uv))^\wedge$, we have
\begin{math}
  \left(\begin{smallmatrix}
    z \\ w
  \end{smallmatrix}\right) =
  R^{-1}_\a
  \left(\begin{smallmatrix}
    f_3 \\ f_4
  \end{smallmatrix} \right)
\end{math}
for some $\a\in]0,\pi[$ and $(f_3,f_4)\in\E_{\l_3}\!\times\!\E_{\l_4}$. If $\a=\frac{\pi}{2}$,
    that is if $(w,z)\in\E_{\l_3}\times\E_{\l_4}$, then signs may be chosen such that
    $X\in\mathfrak{K}'_{1,2,3}(D^{(3)}_k)$. Otherwise, we choose $z$ such that
    $\a<\frac{\pi}{2}$, and thereafter $u,v$ such that
    $X\in\mathfrak{K}'_{1,2}(D^{(3)}_k)$. This gives the second set. \bx

\section{Parametrisation of symmetric matrices} \label{parametrisering}

In this section we parametrise the sets of symmetric matrices having a fixed set of
eigenpairs (Proposition \ref{sym}). In Proposition \ref{flipp} we also give cross-sections
for the orbit sets of these matrices under the $\mathbb{Z}_2$-action realised by changing
the sign of a standard basis vector in the underlying vector space $\R^n$.

Let
\begin{equation} \label{D}
  D=
  \begin{pmatrix}
    \l_1\I_{s_1}\\
      &\ddots\\
      &&\l_k\I_{s_k}
  \end{pmatrix}
\in\R^{n\times n}
\end{equation}
where $\l_i\neq\l_j$ if $i\neq j$.
We write $\p=\p(D)$ for the set of eigenpairs of $D$.

By the real spectral theorem, the map
$\mathfrak{c}:\mathrm{O}_n(\R)\rightarrow \mathrm{Sym}_{\p},\:T\mapsto T^{-1}DT$ is surjective.
In order to get an irredundant parametrisation of $\mathrm{Sym}_{\p}$, one might try to find a
cross-section $\mathcal{S}\subset \mathrm{O}_n(\R)$ for the set
$\mathcal{P}=\{\mathfrak{c}^{-1}(A)\subset \mathrm{O}_n(\R)\mid A\in \mathrm{Sym}_{\p}\}$ of
preimages of $\mathfrak{c}$ in $\mathrm{O}_n(\R)$.

We denote by $C(D)$ the centraliser in $\mathrm{O}_n(\R)$ of the matrix $D$.
For any $S,T\in \mathrm{O}_n(\R)$, we have
\begin{equation}
  S^{-1}DS=T^{-1}DT\Leftrightarrow DST^{-1}=ST^{-1}D\Leftrightarrow 
  \exists N\in C(D):\;S=NT .
\end{equation}
This means that $\mathcal{P}=C(D)\backslash\mathrm{O}_n(\R)$, the set of left cosets of
$C(D)$ in $\mathrm{O}_n(\R)$. 
Note that $C(D)\subset\mathrm{O}_n(\R)$ is in general not a normal subgroup and
therefore, left and right cosets not necessarily coincide.

Let $v\in\R^n$ be a unit vector. There exists a unique $\tau_1\in[0,\pi]$ such that
$v=\cos\tau_{1}e_1+\sin\tau_{1}v_2$ for some $v_2\in e_1^\perp$ of unit length. Continuing
this procedure, in $n$ steps we get the polar coordinates $p(v)=(\tau_1,\ldots,\tau_n)$
for $v$, for which $\tau_n\in\{0,\pi\}$, $v_n=\cos\tau_n e_n=\pm e_n$,
$v_i=\cos\tau_{i}e_i+\sin\tau_{i}v_{i+1}$ and finally $v=v_1$. Adding the condition that
$\tau_i=0$ whenever $\tau_{i-1}\in\{0,\pi\}$, the coordinates of any unit vector
$v\in\R^n$ are unique. We denote by
\begin{equation*}
  \mathfrak{T}_n=\left\{(\tau_1,\ldots,\tau_n)\in[0,\pi]^{n-1}\!\times\!\{0,\pi\} \mid
  \tau_{i-1}\in\{0,\pi\}\Rightarrow\tau_{i}=0\right\}
\end{equation*}
the set of polar coordinates on the unit sphere of $\R^n$.

Given $m\in\underline{n}$ and $\tau\in\mathfrak{T}_m$, set
$y=p^{-1}(\tau_2,\ldots,\tau_m)\in\R^{m-1}$ and define
\begin{align*}
  \tilde{R}(\tau)&= \left(
  \begin{array}{c|c}
    \cos\tau_{1} & -\sin\tau_{1} y^\ast\\
    \hline
    \sin\tau_{1} y & \I_{m-1}-(1-\cos\tau_{1})yy^\ast
  \end{array}
  \right)\in\mathrm{O}_m(\R), \\
  R(\tau)&=
  \begin{pmatrix}
    \I_{n-m} \\ & \tilde{R}(\tau)
  \end{pmatrix}\in\mathrm{O}_n(\R).
\end{align*}
If $\tau_1\not\in\{0,\pi\}$, then $\tilde{R}(\tau)$ is the matrix of rotation in
$\sp\{e_1,Te_1\}\subset\R^m$ mapping $e_1\mapsto p^{-1}(\tau)$. Otherwise,
$\tilde{R}(\tau)=\I_m$ if $\tau_1=0$ and $\tilde{R}(\tau)=\I_m -2e_1e_1^\ast$ if $\tau_1=\pi$.

Fix a matrix $T\in\mathrm{O}_n(\R)$ and let
$\tau^n=(\tau^n_1,\ldots,\tau^n_n)=p(Te_1)$.
Now $Te_1=R(\tau^n)e_1$.
Since $T$ is orthogonal, $(Te_1)^\perp=Te_1^\perp=T(\sp\{e_2,\ldots,e_n\})$. The matrix
$T$ therefore factors uniquely as 
$T=R(\tau^n)\left(\!\begin{smallmatrix}1\\&T_{1}\end{smallmatrix}\!\right)$, where
$T_{1}\in \mathrm{O}_{n-1}(\R)$. Setting $\tau^{n-1}=p(T_1e_1)$ we have
$T_1e_1=\tilde{R}(\tau^{n-1})e_1$ and $Te_2=R(\tau^n)R(\tau^{n-1})e_2$. Proceeding inductively, we
get a factorisation $T=R(\tau^n)\cdots R(\tau^1)$, where $\tau^m\in\mathfrak{T}_m$ for all
$m\in\underline{n}$.

Conclusively,
\begin{equation} \label{bijection}
  \mathfrak{r}:\mathfrak{T}_1\times\cdots\times\mathfrak{T}_n\rightarrow \mathrm{O}_n(\R),\;
  (\tau^1,\ldots,\tau^n)\mapsto R(\tau^n)\cdots R(\tau^1)
\end{equation}
is a bijection. For simplicity, we shall use the notation $R_i=R(\tau^{n+1-i})$ when 
$\tau=(\tau^1,\ldots,\tau^n)$ is given.

We now make a series of important (and, inevitably, very technical) def\-in\-itions. Let
$i,j\in\underline{n}$.
First, if $i<j$, then for any $\tau\in\mathfrak{T}_1\times\cdots\times\mathfrak{T}_n$ we
set 
\begin{equation*}
  \sigma_\tau(i,j)=\sum_{r=1}^{n-j+1}\cos\tau^{n-j+1}_r\cos\tau^{n-i+1}_{j-i+r}
  \prod_{m=1}^{r-1}\sin\tau^{n-j+1}_m\sin\tau^{n-i+1}_{j-i+m}.
\end{equation*}
Second, we define $N(i,j)$ to be the set of all
$\tau\in\mathfrak{T}_1\times\cdots\times\mathfrak{T}_n$ satisfying the following conditions:
\begin{enumerate}
\item If $i\geqslant j$ then either $\tau^{n-j+1}_{i-j+1}=\frac{\pi}{2}$ or 
  $\exists k\in\underline{i-j}:\tau^{n-j+1}_k\in\{0,\pi\}$.
\item If $i<j$ then $\sigma_\tau(i,j)=0$.
\end{enumerate}
Third, we take $P(i,j)$ to be the set of all
$\tau\in\mathfrak{T}_1\times\cdots\times\mathfrak{T}_n$ for
which the following hold true:
\begin{enumerate}
\item If $i\geqslant j$ then $\tau^{n-j+1}_{i-j+1}\in\left[0,\frac{\pi}{2}\right]$.
\item If $i<j$ then $\sigma_\tau(i,j)\leqslant 0$.
\end{enumerate}

In the case $i\geqslant j$, an element
$\tau\in\mathfrak{T}_1\times\cdots\times\mathfrak{T}_n$ belongs to $N(i,j)$ if and only if
$e_i^{\ast}R_je_j=0$, and to $P(i,j)$ if and only if $e_i^{\ast}R_je_j\geqslant0$. 
If $i<j$, we have $(R_ie_i)^{\ast}R_je_j=\sigma_\tau(i,j)\prod_{m=1}^{j-i}\sin\tau^{n-i+1}_{m}$. 
Also, if $R_ie_i\in\sp\{e_i,\ldots,e_{j-1}\}$ then $\sigma_\tau(i,j)=0$, so
$-(R_ie_i)^{\ast}R_je_j\geqslant0$ precisely when $\sigma_\tau(i,j)\leqslant0$, and
$-(R_ie_i)^{\ast}R_je_j=0$ precisely when $\sigma_\tau(i,j)=0$.
Hence $\tau\in N(i,j)$ is equivalent to $-(R_ie_i)^{\ast}R_je_j=0$, and
$\tau\in P(i,j)$ is equivalent to $-(R_ie_i)^{\ast}R_je_j\geqslant0$.

Let $D$ be a matrix of the form given by (\ref{D}), and let $\mathfrak{K}(D)$ be the set
of $\tau\in\mathfrak{T}_1\times\cdots\times\mathfrak{T}_n$ satisfying the following
property: \\
For all $j\in\underline{n-1}$, $r\in\underline{k}$ and all $i\in\mathbb{N}$
  satisfying $\sum_{\mu=1}^{r-1}s_\mu\!<\!i\!<\!\sum_{\nu=1}^{r}s_\nu$, the following
  implications hold true:
\begin{enumerate}
\item If $\tau\in \bigcap_{t=1}^{j-1}N(i,t)$ then $\tau\in P(i,j)$.
\item If $\tau\in \left(\bigcap_{t=1}^{j-1}N(i,t)\right)\minus N(i,j)$, then
   $\tau\not\in N(i-1,h)$ for some
  $h\in\underline{j-1}$. 
\end{enumerate}

\begin{prop}\label{sym}
The restricted map $\mathfrak{c}\circ\mathfrak{r}:\mathfrak{K}(D)\rightarrow \mathrm{Sym}_{\p}$ is a
bijection.
\end{prop}

\textit{Proof:} We need to show that $\mathfrak{r}(\mathfrak{K}(D))$ is a cross-section for
$C(D)\backslash\mathrm{O}_n(\R)$. The centraliser $C(D)$ of $D$ consists of those orthogonal
matrices that leave all eigen\-spaces of $D$ invariant.

Let $T\in \mathrm{O}_n(\R)$. Consider $(\l,s)\in{\p}$ with
$\E_\l(D)=\sp\{e_l,\ldots,e_{l+s-1}\}$.
We denote by $v_{i}^\l$ the vector in $\R^s$ given by the restriction of the $i$th column of
$T$ to the rows $l$ to $l+s-1$ inclusive, \ie, $v_{i}^\l=P_{\E_\l(D)}T(e_i)$. 
For $j\in\underline{s}$ we define 
$\imath(j)=\min\{i\mid \dim(\sp\{v_1^\l,\ldots,v_i^\l\})=j\}$.
Since $T$ is invertible, $\dim(\sp\{v_1^\l,\ldots,v_n^\l\})=s$ and the definition of
$\imath:\underline{s}\to\underline{n}$ is consistent.
Now let $S_{\l}\in \mathrm{O}_s(\R)$ be the matrix defined by
$e_1=(S_{\l}v^\l_{\imath(1)})^\wedge$ and
$e_{j+1}=
\left(S_{\l}P_{\{e_1,\ldots,e_{\imath(j)}\}^\perp}(v^\l_{\imath(j+1)})\right)^\wedge$.

Let
\begin{equation*}
  S_T=\left(\!\begin{smallmatrix}S_{\l_1}\\&\ddots\\&&S_{\l_k}\end{smallmatrix}\!\right).
\end{equation*}
Clearly, $S_T\in C(D)$.
The matrix $P_{\E_\l(D)}S_TT$ (that is the restriction of $S_TT$ to the rows $l$ to
$l+s-1$) now has the following form:
\begin{equation}\label{normalmatris}
  P_{\E_\l(D)}S_TT=
  \left(\!
  \begin{array}{cccccccccccccccccc}
    0&\cdots&0& \bullet &\ast&\cdots \\
    0&\cdots&&&&\cdots&0& \bullet &\ast&\cdots \\
    \vdots & \ddots \\
    &&&&&&&&&&&&&\ddots\\
    0&\cdots&&&&&&&&\cdots&0&\bullet&\ast&\cdots
  \end{array}
  \!\right).
\end{equation}
Here a bullet indicates a strictly positive entry, and an asterisk denotes an
arbitrary real number. The first non-zero entry of each row will be positive, and
located strictly to the right of the first non-zero entry of the previous row. These
elements we call pivotal elements. The matrix
$S_TT$ is a tower of blocks of this type, one for each eigenvalue $\l$ of $D$. Clearly,
$S_TT$ is an exact invariant for the left cosets of $C(D)$ to which it belongs. 

If we read $\tau\in N(i,j)$ as
``the entry $(i,j)$ of the matrix $\mathfrak{r}(\tau)$ is zero''
and $\tau\in P(i,j)$ as
``the entry $(i,j)$ of $\mathfrak{r}(\tau)$ is non-negative''
we see that $\mathfrak{r}(\mathfrak{K}(D))$ is precisely the set of all orthogonal
matrices of the above form. Although these descriptions are not true in general, we will
show that they are correct if all entries to the left of $(i,j)$ in the $i$th row are zero,
which is the case in the definition of $\mathfrak{K}(D)$.

Let $\tau\in\mathfrak{T}_1\times\cdots\times\mathfrak{T}_n$.
Fix $i,j\in\underline{n}$, and assume $\tau\in\bigcap_{t=1}^{j-1}N(i,t)$.
The matrix $T=\mathfrak{r}(\tau)$ factors into a product of rotations and reflections
$T=R_1\cdots R_n=R(\tau^n)\cdots R(\tau^1)$. Since
$R_m=\left(\!\begin{smallmatrix}\I_{m-1}\\&\tilde{R}_m\end{smallmatrix}\!\right)$ we have
\begin{math}
  e_i^{\ast}Te_j=e_i^{\ast}R_1\cdots R_ne_j=e_i^{\ast}R_1\cdots R_je_j.
\end{math}

On the other hand, for all $t<\min\{i,j\}$ we have
\begin{equation} \label{pilar}
\tau\in N(i,t)\;\stackrel{i\geqslant t}{\Rightarrow}\;e_i^{\ast}R_te_t=0
\;\stackrel{t\neq i}{\Rightarrow}\;R_te_i=e_i\;\Rightarrow\;R_t^{\ast}e_i=e_i
\;\Rightarrow\;e_i^{\ast}R_t=e_i^{\ast}\,.
\end{equation}
Hence $e_i^{\ast}Te_j=e_i^{\ast}R_m\cdots R_je_j$, where $m=\min\{i,j\}$.

In the case $i\geqslant j$, this means that $e_i^{\ast}Te_j=e_i^{\ast}R_je_j$. Therefore
the statement $\tau\in P(i,j)$ is equivalent to $e_i^{\ast}Te_j\geqslant0$ and 
$\tau\in N(i,j)$ is equivalent to $e_i^{\ast}Te_j=0$, under the present condition that
$\tau\in\bigcap_{t=1}^{j-1}N(i,t)$. 

As for the case $i<j$, we get $e_i^{\ast}Te_j=e_i^{\ast}R_i\cdots R_je_j$. Since by
assumption $e_i^{\ast}R_ie_i=0$, $R_i$ must be a rotation, and it follows that
$R_i^{\ast}e_i=-R_ie_i$. Moreover, $\tau\in N(i,t)$ for all $t\in\{i\!+\!1,\ldots,j-1\}$
implies that $(R_ie_i)^{\ast} R_t=(R_ie_i)^{\ast}$, by arguments analogous to (\ref{pilar}).
We calculate
\begin{equation*}
  e_i^{\ast}Te_j=e_iR_i\cdots R_je_j=-(R_ie_i)^{\ast}R_{i+1}\cdots R_je_j=-(R_ie_i)^{\ast} R_je_j\,.
\end{equation*}
Hence $\tau\in N(i,j)$ if and only if $e_i^{\ast}Te_j=0$ and $\tau\in P(i,j)$ if and
only if $e_i^{\ast}Te_j\geqslant0$. 

Summarising the above, we get that $\mathfrak{K}(D)$
parametrises the set of orthogonal matrices for which the blocks $P_{\E_\l(D)}T$ have
the form $(\ref{normalmatris})$. This proves the proposition. \bx

Let $i\in\underline{n}$. We write $\mathbb{Z}^{(i)}_2=\{\I_n,\Sigma_i\}$, where
$\Sigma_i=\I_n-2e_ie_i^{\ast}$ is the matrix of reflection in the hyperplane $e_i^\perp$. The
(multiplicative) group $\mathbb{Z}^{(i)}_2$ acts on $\mathrm{Sym}_\p$ by $g\cdot B=gBg$.

For any $m\in\underline{n+1-i}$, set $\check{m}=i-1+m$. Given $\tau\in\mathfrak{K}(D)$, we define
\begin{align*}
  L_i(\tau)&= \left\{ m\in\underline{n+1-i} \;\Big|\; \tau^{n+1-i}_m\neq\frac{\pi}{2} 
  \mbox{ and } \left(\, \tau^{n+1-i}_{m-1}\not\in\{0,\pi\} \mbox{ if } m>1 \,\right) \right\}\\ 
  &= \left\{m\in\underline{n+1-i} \;\big|\; e^\ast_m \,p^{-1}(\tau^{n+1-i})\neq0 \right\} \\ 
\intertext{and}
K_i(\tau)&= \left\{m\in L_i(\tau) \;\bigg|\; \tau\in \bigcap_{t=1}^{i-1}N(\check{m},t) \right\}
\subset L_i(\tau).
\end{align*}
The set $\{(\check{m},i)\}_{m\in K_i(\tau)}$ is the set of all pivotal elements on and
below the diagonal in the $i$th column of $\mathfrak{r}(\tau)$.

Let $\mathfrak{K}_i(D)\subset\mathfrak{K}(D)$ be the set of all $\tau\in\mathfrak{K}(D)$
such that either
\begin{enumerate}
\item $L_i(\tau)\minus K_i(\tau)\neq\emptyset$ and $\tau^{n+1-i}_l<\frac{\pi}{2}$ where
  $l=\min(L_i(\tau)\minus K_i(\tau))$, or
\item $L_i(\tau)=K_i(\tau)$ and $|L_i(\tau)|=1$, or
\item $L_i(\tau)=K_i(\tau)$, $|L_i(\tau)|>1$ and $\tau\in P(\check{r},t)$ where 
  $r=\min L_i(\tau)$ and
  $t=\min \left\{ m\in\underline{n} \mid m>i,\; \tau\not\in N(\check{r},m) \right\}$.
\end{enumerate}

\begin{prop} \label{flipp}
  For any $i\in\underline{n}$, the set
  $(\mathfrak{c}\circ\mathfrak{r})(\mathfrak{K}_i(D))$ is a cross-section for
  $\mathrm{Sym}_{\p}\slash\mathbb{Z}^{(i)}_2$.
\end{prop}

\textit{Proof:}
Consider $B=T^{-1}DT\in\mathrm{Sym}_\p$ with $T=\mathfrak{r}(\tau)$ and
$\tau\in\mathfrak{K}(D)$. 
The element $\Sigma_i$ acts as 
$\Sigma_i\cdot B=\Sigma_iB\Sigma_i=\Sigma_iT^{-1}DT\Sigma_i=(T\Sigma_i)^{-1}D(T\Sigma_i)$.
Suppose $T\Sigma_i=\mathfrak{r}(\tilde{\tau})$ and
$S_{T\Sigma_i}T\Sigma_i=\mathfrak{r}(\bar{\tau})$, where
$S_{T\Sigma_i}\in C(D)$ is the unique matrix, given in the proof of
Proposition~\ref{sym}, for which $S_{T\Sigma_i}T\Sigma_i\in\mathfrak{K}(D)$.
We need to show that for every $\tau\in\mathfrak{K}(D)$ either $\tau\in\mathfrak{K}_i(D)$
or $\bar{\tau}\in\mathfrak{K}_i(D)$, but not both.

Certainly, $T\Sigma_ie_i=-Te_i$ and $T\Sigma_ie_j=Te_j$ for all $j\neq i$.
From this follows that $S_{T\Sigma_i}=\prod_{h\in H}\Sigma_h$, where $H$ is
the set of all $h\in\underline{n}$ for which the element $(h,i)$ of $T$ is a pivotal element.

On the other hand,
\begin{align*}
  &T\Sigma_i=R(\tau^n)\cdots R(\tau^1)\Sigma_i=R(\tau^n)\cdots R(\tau^{n+1-i})\Sigma_i
  R(\tau^{n-i})\cdots R(\tau^1), \\
  &R(\tau^{n+1-i})\Sigma_i=
  \begin{pmatrix}
    \I_{i-1}\\
    &\tilde{R}(\tau^{n+1-i})
  \end{pmatrix}
  \Sigma_i=
  \begin{pmatrix}
    \I_{i-1}\\
    &\tilde{R}(\tau^{n+1-i})\Sigma_1
  \end{pmatrix}.
\end{align*}
Setting $\theta=\tau^{n+1-i}_{1}$ and
$y=p^{-1}(\tau^{n+1-i}_2,\ldots,\tau^{n+1-i}_{n+1-i})$, we get
\begin{multline*}
  \tilde{R}(\tau^{n+1-i})\Sigma_1= \\ =
  \begin{pmatrix}
    \cos\theta & -\sin\theta y^{\ast} \\
    \sin\theta y  & \I_{n-1}-(1-\cos\theta)yy^{\ast}
  \end{pmatrix} \Sigma_1 =
  \begin{pmatrix}
    -\cos\theta & -\sin\theta y^{\ast} \\
    -\sin\theta y  & \I_{n-1}-(1-\cos\theta)yy^{\ast}
  \end{pmatrix} = \\ 
  \rule{0cm}{6mm} =
  \begin{pmatrix}
    \cos(\pi-\theta) & -\sin(\pi-\theta)(-y)^{\ast} \\
    \sin(\pi-\theta)(-y) & \I_{n-1}-(1-\cos(\pi-\theta))(-y)(-y)^{\ast}
  \end{pmatrix} A_y
  =  R(\tilde{\tau}^{n+1-i})A_y
\end{multline*}
where
\begin{equation*}
  A_y=\begin{pmatrix}1\\&\I_{n-1}-2yy^{\ast}\end{pmatrix} \mbox{ and }
  \tilde{\tau}^{n+1-i}=
  p\begin{pmatrix}\cos(\pi-\theta)\\\sin(\pi-\theta)(-y)\end{pmatrix} .
\end{equation*}
Hence $\tilde{\tau}^{n+1-i}_l=\pi-\tau^{n+1-i}_l$ for all $l$ for which
$\tau^{n+1-i}_{l-1}\not\in\{0,\pi\}$.
We now get three cases, as in the definition of $\mathfrak{K}_i(D)$.

If $L_i(\tau)\minus K_i(\tau)\neq\emptyset$: Let $l=\min(L_i(\tau)\minus K_i(\tau))$.
Now $\tilde{\tau}_l^{n+1-i}=\pi-\tau_l^{n+1-i}$. Since $l\not\in K_i(\tau)$ we have
$\check{l}\not\in H$ and therefore $\bar{\tau}_l^{n+1-i}=\tilde{\tau}_l^{n+1-i}$. 
Hence precisely one of $\tau_l^{n+1-i}$ and $\bar{\tau}_l^{n+1-i}$ belongs to the
interval $\left[0,\frac{\pi}{2}\right[$. This is the first case in the definition of
$\mathfrak{K}_i(D)$. 

If $L_i(\tau)\minus K_i(\tau)=\emptyset$:
Then the $i$th column of $T$ equals $p^{-1}(\tau^{n+1-i})$ and
every non-zero entry is a pivotal element. 
If $|L_i(\tau)|=1$, then $p^{-1}(\tau^{n+1-i})=(\delta_{hj})_{j\in\underline{n}}$ for some
$h\in\underline{n}$. Hence the $h$th row of $T$ looks like
$(\delta_{ij})_{j\in\underline{n}}^{\ast}$. This means that
$S_{T\Sigma_i}T\Sigma_i=\Sigma_hT\Sigma_i=T$ and consequently $\bar{\tau}=\tau$.
This is the second case.

Otherwise, if $|L_i(\tau)|=|K_i(\tau)|>1$, set $r=\min L_i(\tau)$. 
The element $(\check{r},i)$ of $T$ is not equal to $\pm1$, so the $\check{r}$th row must contain some
other non-zero element. Say the first such element is $(\check{r},t)$. 
Note that then $t>i$, since $(\check{r},i)$ is a pivotal element of $T$.
We have $S_{T\Sigma_i}T\Sigma_i=\left(\prod_{l\in L_i(\tau)}\Sigma_{\check{l}}\right)T\Sigma_i$. 
Since $r\in L_i(\tau)$, left multiplication with $S_{T\Sigma_i}$ will change the sign of the
element $(\check{r},t)$, whereas right multiplication with $\Sigma_i$ leaves it unchanged. This
implies that precisely one of $\tau$ and $\bar{\tau}$ belongs to $P(r,t)$, which is
the requirement in the third case of the definition of $\mathfrak{K}_i(D)$.  \bx

Finally, we remark that $\mathfrak{K}_{i_1,\ldots,i_l}(D)=
\bigcap_{j=1}^l\mathfrak{K}_{i_j}(D)$ is indeed a cross-section for the action of the direct
product $\prod_{j=1}^l\mathbb{Z}^{(i_j)}_2$ on $\mathrm{Sym}_{\p}$ by
$(g_1,\ldots,g_l)\cdot B=g_1\cdots g_lBg_l\cdots g_1$.

\bibliographystyle{plain}
\bibliography{../litt.bib}

\end{document}